\documentclass[a4paper,12pt]{article} 

\usepackage{amsmath,amssymb,fullpage,indentfirst}
\usepackage{enumerate}
\usepackage{graphicx}
\usepackage{subfigure}
\usepackage{hyperref}
\usepackage{url}
\usepackage[utf8]{inputenc}
\usepackage[T1]{fontenc}

\DeclareMathOperator{\arccot}{arccot}

\newtheorem{theorem}{Theorem}
\renewcommand{\textit}[1]{\textsl{#1}}
\title{Revisiting Gilbert Strang's ``A Chaotic Search for $i$''}
\author{Ao Li and Robert M. Corless}
\date{}

\begin{document}
	
	\maketitle
	
	\begin{abstract}
		In the paper ``A Chaotic Search for $i$''~(\cite{strang1991chaotic}), Strang completely explained the behaviour of Newton's method when using real initial guesses on $f(x) = x^{2}+1$, which has only a pair of complex roots $\pm i$. He explored an exact symbolic formula for the iteration, namely $x_{n}=\cot{ \left( 2^{n} \theta_{0} \right) }$, which is valid in exact arithmetic. In this paper, we extend this to to $k^{th}$ order Householder methods, which include Halley's method, and to the secant method. Two formulae, $x_{n}=\cot{ \left( \theta_{n-1}+\theta_{n-2} \right) }$ with $\theta_{n-1}=\arccot{\left(x_{n-1}\right)}$ and $\theta_{n-2}=\arccot{\left(x_{n-2}\right)}$, and $x_{n}=\cot{ \left( (k+1)^{n} \theta_{0} \right) }$ with $\theta_{0} = \arccot(x_{0})$, are provided. The asymptotic behaviour and periodic character are illustrated by experimental computation. We show that other methods (Schr\"{o}der iterations of the first kind) are generally not so simple. We also explain an old method that can be used to allow Maple's \textsl{Fractals[Newton]} package to visualize general one-step iterations by disguising them as Newton iterations.
	\end{abstract}
	

		\section{Introduction} 
		
		The study of discrete dynamical systems, denoted generically here by $x_{n+1} = F(x_n)$ with $x_0 \in \mathbb{C}^d$ a $d$-dimensional complex vector and $F$ being a typically nonlinear map, is both old and important in mathematics and its applications. One extremely well-studied aspect of this is the use of such iterations to search for fixed points of the map; if the map is itself of the form $x - D^{-1}(x)\,G(x)$, then if $D$ is not singular at the fixed point, we will have found a zero of the (usually nonlinear) map $G(x)$. Finding zeros and equilibria is of course an important question in many applications, such as design or game theory.
		
		It may seem surprising that the study of just the simplest nonlinear example---even just in one dimension---namely $f(x) = x^2 +1$ and various iteration schemes to solve it, such as Newton's method and variations, can clarify deep questions for the general case, but indeed this is so.  For an earlier instance of this, using ideas of Charles M.~Patton and also citing Strang's paper, see~\cite{corless1998newton}.
		
		This paper reports on what began as a student project in a graduate course, Open Problems in Experimental Mathematics; namely trying to extend the results of~\cite{strang1991chaotic} to other iteration methods. After solving the problem, we found the paper~\cite{rhouma2005fibonacci} which had extended the results at least to Halley's method and to the secant method; thus the problem was not as open as we had thought.  However, the extension to \textsl{all} Householder methods, our theorem~\ref{thm:House}, is new to this current paper.
		
		For completeness, this current paper also includes our  rediscovery of the extension of Strang's results to Halley iteration and secant iteration. We then give our main theorem, which extends the results (using a symbolic $n$th derivative) to Householder methods. We then use Maple's \textsl{Fractals} package to show why we believe that Schroeder's first methods are more difficult to understand and likely cannot be explained with a similar trick.  
		
		We begin with a review of Newton's method for finding zeros of $f(x)$.
		
		\section{A review of Newton's Method}
		
		Newton's method and its variants are workhorses of scientific computing: they replace the task of solving $f(x) = 0$ with an iteration $x_{n+1} = F (x_n)$ which maps a ``starting guess'' $x_0$ to a sequence $x_1, x_2, x_3, \dots$ which hopefully quickly converges to a solution $x^*$ such that $f(x^*) = 0$. The basic idea was indeed used by Newton himself, though in a careful context of repeatedly shifting the point of expansion of a finite Taylor series for a polynomial until the first term, $f(x_n)$, became negligibly small. It was Euler who first gave us Newton's method for scalar $f(x)$. [Wanner (\cite{butcher1996runge} and \cite{hairer2003geometric}) tells us that, symmetrically enough,  it was Newton who first used what is now known as the symplectic Euler method. See~\cite{kalantari2008polynomial} for more historical details.]
		
		Schr\"{o}der extended this to all higher orders; his discoveries are continually reinvented (\cite{schroder1998infinitely}), which just seems to be a fact of life even in a modern age where information is easy to find. We will use Schr\"{o}der's point of view to explain Newton's method, below.
		
		Consider
		\begin{equation}
		\begin{aligned}
		y &= f(x_{0} + \varepsilon) \\
		&= f(x_{0}) + f'(x_{0}) \varepsilon + \frac{1}{2} f''(x_{0}) \varepsilon^{2} + \cdots,
		\end{aligned}
		\end{equation}
		assuming $f(x)$ sufficiently differentiable. We now reverse the series, which we can do provided $f'(x_{0}) \ne 0$:
		\begin{equation}
		\varepsilon = \frac{1}{f'(x_{0})}(y - f(x_{0})) + A_{2}(y - f(x_{0}))^{2} + \cdots.
		\end{equation}
		The coefficient $A_{2} = -f''(x_{0}) / (2f'(x_{0})^{3})$ is known in terms of $f$ and its derivatives at $x_{0}$. Formulas are known and tabulated for the first few $A_{k}$, in fact; and effective means are available for computing as many $A_{k}$ as one could desire, although the cost of such computation increases as the desired number of $A_{k}$ increases. This was known already to Lagrange, and one theoretically useful method for finding the $A_{k}$ is called the Lagrange Inversion Formula (\cite{comtet2012advanced}).
		
		To find $x^*=x_{0}+\varepsilon$ such that $y=0$, simply put $y=0$ in the series for $\varepsilon$. If we have all terms, and the series converges, then adding the result $\varepsilon$ to the known $x_{0}$ gives the desired $x^*$.
		
		In practice one truncates the series. For Newton's method, we ignore $A_{2}$ and all subsequent terms and take
		\begin{equation}
		\hat{\varepsilon} = \frac{1}{f'(x_{0})}(0-f(x_{0})),
		\end{equation}
		giving a new estimate $x_{1}=x_{0}+\hat{\varepsilon}$ or 
		\begin{equation}
		x_{1}=x_{0}-\frac{f(x_{0})}{f'(x_{0})}.
		\end{equation}
		Newton's idea is to use this formula repeatedly:
		\begin{equation}
		\begin{aligned}
		x_{2} &= x_{1}-\frac{f(x_{1})}{f'(x_{1})}, \\
		x_{3} &= x_{2}-\frac{f(x_{2})}{f'(x_{2})},
		\end{aligned}
		\end{equation}
		which requires repeated (usually costly) evaluation of $f$ and its derivatives, and comes with no true \textsl{a priori} guarantee of success. Better alternatives are continually sought.
		
		Because the series for Newton's method has an error $O(\varepsilon^{2})$, iterating it will (in the best case) square the previous error, which is called ``quadratic convergence''. 
		
		What if we also keep the $A_{2}$ term? Then,
		\begin{equation}
		\Tilde{\varepsilon} = \frac{1}{f'(x_{0})} (0-f(x_{0})) - \frac{f''(x_{0})}{2f'(x_{0})^{3}} (0-f(x_{0}))^{2}.
		\end{equation}
		So now,
		\begin{equation}
		\begin{aligned}
		x_{1} &= x_{0}+\Tilde{\varepsilon} \\
		&= x_{0}-\frac{f(x_{0})}{f'(x_{0})} - \frac{f''(x_{0})}{2f'(x_{0})^{3}} f^{2}(x_{0}),
		\end{aligned}
		\end{equation}
		and this method is cubically convergent. It has the disadvantage of needing the prior computation of the second derivative $f''$; nonetheless the method is viable.
		
		However, this method is not often used. Instead, another cubically convergent method, known as Halley's method, is used:
		\begin{equation} \label{eqn::Halley}
		x_{n+1} = x_{n} - \dfrac{f(x_{n})}{f'(x_{n}) - \dfrac{f(x_{n})f''(x_{n})}{2f'(x_{n})}}.
		\end{equation}
		If $f(x_{n})$ is small, then
		\begin{equation}
		\begin{aligned}
		\dfrac{1}{f'(x_{n}) - \dfrac{f(x_{n})f''(x_{n})}{2f'(x_{n})}} &= \dfrac{1}{f'(x_{n})} \cdot \dfrac{1}{1 - \dfrac{f(x_{n})f''(x_{n})}{2f'(x_{n})^{2}}} \\
		&= \dfrac{1}{f'(x_{n})} \cdot \left( 1 + \dfrac{f(x_{n})f''(x_{n})}{2f'(x_{n})^{2}} \right) + O(f^{2}(x_{n})),
		\end{aligned}
		\end{equation}
		and we recover the cubic Schr\"{o}der iteration to the same order of error.
		
		Higher-order Schr\"{o}der iterations---indeed methods of arbitrary order---are possible and occasionally useful.
		
		But in fact lower-order methods such as the secant method discussed below, and their multidimensional analogues such as the BFGS method, are cheaper in practice (once they get started) because they re-use more than just the previous iterate (see~\cite{neumaier2001introduction} for a detailed analysis). The secant method uses the iteration:
		\begin{equation} \label{eqn::secant}
		x_{n+1} = x_{n} - \dfrac{f(x_{n})(x_{n}-x_{n-1})}{f(x_{n}) - f(x_{n-1})}.
		\end{equation}
		Here $f'(x_{n})$ has been replaced with the secant approximation.
		
		Other, more sophisticated schemes such as Inverse Quadratic Interpolation (also called the Dekker-Brent algorithm) can be even more effective; see the documentation for Matlab's fzero command. There the idea is to fit a quadratic in $y$ to three iterates $(x_{0}, y_{0})$, $(x_{1}, y_{1})$ and $(x_{2}, y_{2})$ and set $y=0$ in the result; the formulas are complicated to human eyes but effective computationally (when they do not run into trouble). The following formula is taken from~\cite{corless2013graduate}:
		\begin{equation}
		x_{3} = x_{2} + \dfrac{y_{1}y_{2}(x_{0}-x_{2})}{(y_{0} - y_{1})(y_{0} - y_{2})} + \dfrac{y_{0}y_{2}(x_{1}-x_{2})}{(y_{1} - y_{0})(y_{1} - y_{2})}.
		\end{equation}
		
		\section{Failure of Newton's method}
		Although Newton's method is a crucial algorithm in root finding, it has several known flaws.  It can only find one root at a time, and it does not indicate that all roots are found, or that there are no roots.  Indeed, it runs into trouble even for the simplest nonlinear scalar equation, 
		\begin{equation} \label{eqn::equation}
		f(x)=x^{2}+1=0\>,
		\end{equation}
		which has two roots $x=\pm i$. Newton's method gives the recursive equation:
		\begin{equation} \label{eqn::Newton}
		x_{n+1}=x_{n} - \dfrac{x_{n}^{2}+1}{2x_{n}} = \dfrac{1}{2} \left(x_{n}-\dfrac{1}{x_{n}}\right).
		\end{equation}
		Evidently, any sequences generated by~(\ref{eqn::Newton}) that start from a real number cannot converge to either of the complex points $x=\pm i$, because the iterates must remain real. 
		
		Strang studied these sequences in~\cite{strang1991chaotic}. He recognized a trigonometric identity which is similar to the recursive formula~(\ref{eqn::Newton}), namely
		\begin{equation} 
		\cot (2\theta) = \dfrac{1}{2} \left(\cot \theta -\dfrac{1}{\cot \theta}\right).
		\end{equation}
		If $x_{n}$ is the cotangent of an angle $\theta_{n}$, then the next step gives the cotangent of the double angle $2\theta_{n}$. Therefore, one analytical expression for $x_{n}$ provided by Strang is,
		\begin{equation} \label{eqnN}
		x_{n}=\cot \left( 2^{n} \theta_{0} \right), \quad \text{given} \quad x_{0}=\cot \theta_{0}.
		\end{equation}
		Since $-\infty < \cot \theta < \infty$ for $0 < \theta < \pi$, for any real initial guess, one can uniquely choose $\theta_{0} \in (0,\pi)$ and then analyse the asymptotic behaviour of the sequence. The following results are given in~\cite{strang1991chaotic}. Notice that we may take $2^{n} \theta_{0}$ modulo $\pi$ because $\cot(\phi + k\pi) = \cot \phi$ for $k \in \mathbb{Z}$. 
		\begin{enumerate}[(i)]
			\item If $\theta_{0}=\dfrac{k \pi}{2^{n}}$ for some $k \in \mathbb{Z}$, such as $\theta_{0}=\dfrac{\pi}{4}$, then $x_{n}=\cot{\left( k\pi \right)}$. The iteration blows up, because cotangent is singular at multiples of $\pi$.
			\item If $\theta_{0}=\dfrac{p}{q} \pi$ for any fraction $\dfrac{p}{q}$ other than $\dfrac{k}{2^{n}} \left(k \in \mathbb{Z}\right)$, then the iteration eventually cycles. In addition, when $\theta_{0}=\dfrac{k \pi}{2^{n}-1}$ for some $k \in \mathbb{Z}$, such as $\theta_{0}=\dfrac{\pi}{3}$, we will see $x_{n}=x_{0}$. The iteration of period $n$ cycles from the start point.
			\item If $\theta_{0}=b \pi$ for some irrational number $b$, the iteration is not periodic (or convergent).
		\end{enumerate} 
		
		The map $\theta \to 2\theta$ mod $\pi$ is a variation of the Bernoulli shift map, and well known to be chaotic~(\cite{billingsley1965ergodic}). 
		
		\subsection{The effect of floating-point}
\begin{figure}[h]
			\hfil \includegraphics[width=14cm,height=6cm]{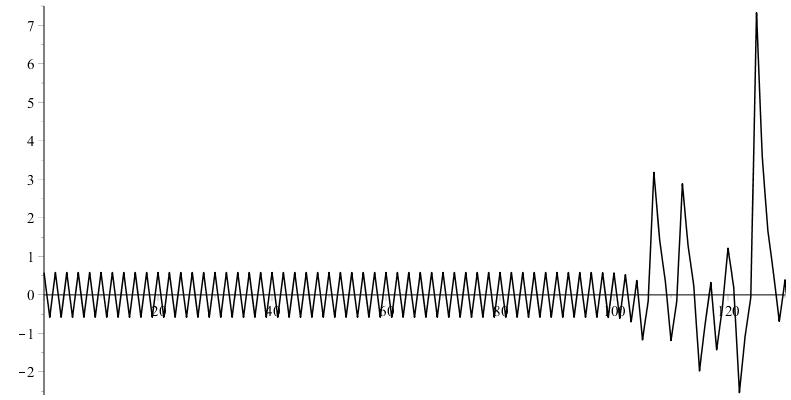}
			\caption{Choose $\theta_{0}= \dfrac{\pi}{3} =\dfrac{\pi}{2^{2}-1}$. The prime period of the sequence is $2$. Numerically, the oscillation is destroyed by the growing round-off error. This happens no matter how many digits are used, although it takes more iterations before the periodicity is destroyed if more digits are used. }
			\label{Ncycle}
		\end{figure}

		Figure~\ref{Ncycle} shows a periodic iteration starting from $x_{0}=\cot (\pi/3)=\sqrt{3}/3$. By simple computation, we know that the sequence oscillates between $\sqrt{3}/3$ and $-\sqrt{3}/3$. However, round-off error interferes if we use floating-point arithmetic. Using Maple and keeping $32$ digits, the periodicity is eventually destroyed by the growing round-off error. 
		
\begin{figure}[h]
			\hfil 
			\includegraphics[width=14cm,height=6cm]{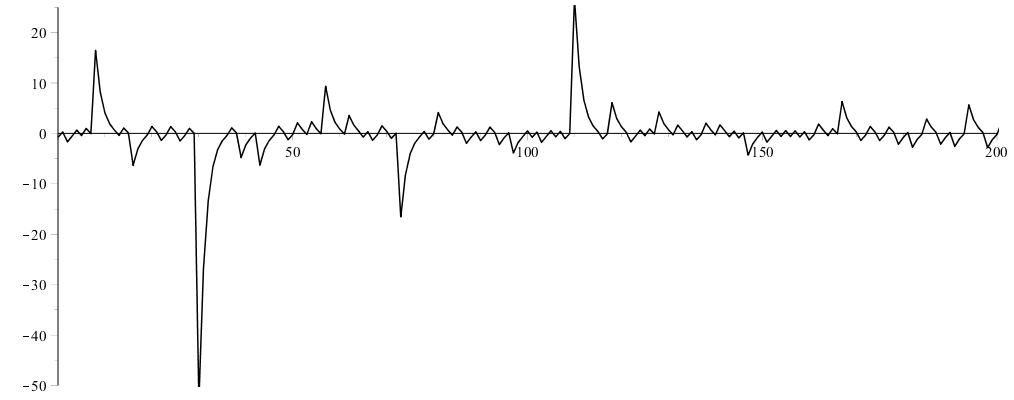}
			\caption{Choose $\theta_{0}=\dfrac{\sqrt{2}\pi}{2}$. The iteration is not periodic.}
			\label{Nchaos}
		\end{figure}

		Figure~\ref{Nchaos} gives an aperiodic example with the initial angle $\theta_{0}=\sqrt{2}\pi/2$. This erratic behaviour in floating-point is not surprising, because the map is chaotic. 
		In the next section, we still use the example $f(x)=x^{2}+1=0$ and extend the result to other algorithms with different orders of convergence, such as Halley's method, the secant method and general Householder's method.
		
		\section{Other root-finding methods }	
		
		\subsection{Halley's method} 
		
		By direct computation, the first and second derivatives of function $f(x)=x^{2}+1$ are $f'(x)=2x$ and $f''(x)=2$. Substituting into iteration~(\ref{eqn::Halley}) gives
		\begin{equation}
		x_{n+1} = \dfrac{x_{n}^{3}-3x_{n}}{3x_{n}^{2}-1}.
		\end{equation}
		
		Inspired by Strang's idea, we also try the trigonometric identities for a match. The formulae we find are
		\begin{equation}
		\cot{(3\theta)}=\dfrac{\cot^{3}{\theta}-3\cot{\theta}}{3\cot^{2}{\theta}-1} \quad \text{and} \quad \tan{(3\theta)}=\dfrac{\tan^{3}{\theta}-3\tan{\theta}}{3\tan^{2}{\theta}-1}.
		\end{equation}
		
		To be similar with the formula found by Strang, we use the cotangent one. Since $\tan(\dfrac{\pi}{2} - \theta) = \cot \theta$, this merely amounts to relabelling the angles. Hence, if $x_{n}=\cot{\theta_{n}}$, then $x_{n+1}=\cot{(3\theta_{n})}$. Then, 
		\begin{equation} \label{eqnH}
		x_{n}=\cot{\left(3^{n}\theta_{0}\right)}, \quad \text{given} \quad x_{0}=\cot \theta_{0}.
		\end{equation} 
		The angle grows exponentially. Compare this formula and expression~(\ref{eqnN}) found by Strang. The only difference is the constants, $2$ for Newton's method and $3$ for Halley's method. This is interesting because it is well-known that the iterates converge quadratically and cubically, respectively, when they converge.
		
		Since formula~(\ref{eqnH}) is close to that for Newton's method, it is natural to see that the iteration displays similar behaviour. 
		
		\paragraph{Case 1} The iteration diverges to infinity. Given $\theta_{0}= \dfrac{k\pi}{3^{n}} (k \in \mathbb{N})$, we see that $\theta_{n}=3^{n}\theta_{0}$ is a multiple of $\pi$, whose cotangent is infinite. Take $\theta_{0}= \pi/9$. Then $x_{1}=\cot(\pi/3)= \sqrt{3}/3$ and $x_{2}=\cot \pi = -\infty$. Doing the iteration numerically, the impact of round-off error arises, leading to a totally different pattern of the sequence. Instead of returning negative infinity as expected, the iteration after two steps gives a very large number $x_{2}=2.1994295969128600552729477352456\times 10^{31}$ by Maple when keeping $32$ digits. The result can be much larger if we use more digits. We also notice that $x_{3}$ is close to one-third of $x_{2}$ and $x_{4}$ is close to one-third of $x_{3}$. This is because that $x_{n+1} = \dfrac{x_{n}^{3}-3x_{n}}{3x_{n}^{2}-1} \approx \dfrac{x_{n}}{3}$ when $x_{n}$ is large. 
		
		\paragraph{Case 2} If exact arithmetic is used, the iteration eventually cycles if $\theta_{0}=\dfrac{p}{q} \pi$ for any fraction $\dfrac{p}{q}$ other than those we mentioned in case 1. Given this initial point, $\theta_{n}=3^{n}p \cdot \dfrac{\pi}{q}$. The sequence oscillates because the denominator $q$ remains as is and $3^{n}p$ modulo $q$ are bounded.
		
		
		\paragraph{Case 2a} The iteration cycles from the start point. To find period-$n$ cycles, require that $x_{n}=x_{0}$, which yields $\theta_{n} \equiv \theta_{0} \mod \pi$. Thus $\theta_{0}= \dfrac{k\pi}{3^{n}-1}$ for some $k \in \mathbb{Z}$.
		
		\begin{figure}[h]
			\centering
			\begin{minipage}[b]{14cm}
				\includegraphics[width=14cm,height=5.6cm]{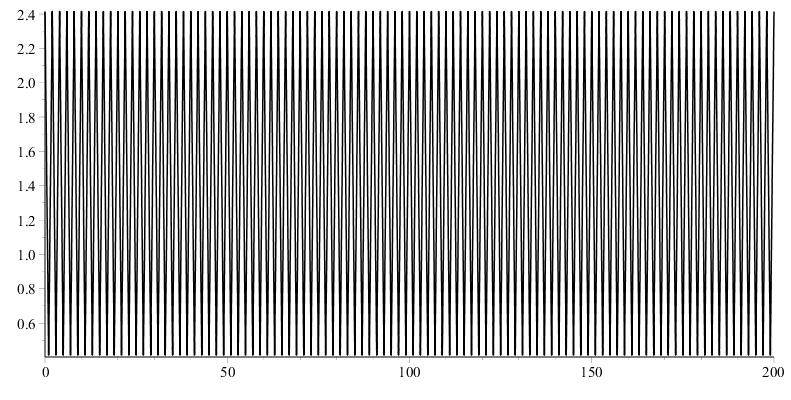} 
				\caption{Period-$2$ cycle given $\theta_{0}= \dfrac{\pi}{8}$. Roundoff error does not seem to bother this instance.} 
				\label{Hcycle1}
				\includegraphics[width=14cm,height=5.6cm]{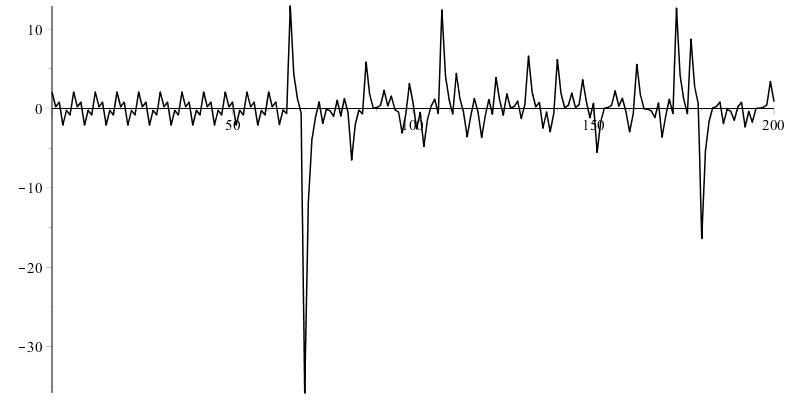}
				\caption{Period-$6$ cycle given $\theta_{0}= \dfrac{\pi}{7}$. Roundoff error quickly destroys the actual periodicity.}
				\label{Hcycle2}
			\end{minipage}
		\end{figure}
		
		Figures~\ref{Hcycle1} and~\ref{Hcycle2} show two examples. When $\theta_{0}= \pi/8$, it is a period-$2$ cycle because $\dfrac{1}{8}= \dfrac{1}{3^{2}-1}$. When $\theta_{0}= \pi/7$, it is a period-$6$ cycle because $\dfrac{1}{7}= \dfrac{104}{3^{6}-1}$. The numerical results are different. The first one looks fine at least for the first $200$ steps. It oscillates exactly between the two limits. However, the second sequence only keeps its periodicity for no more than ten periods (around $60$ steps) before destroyed by the growing round-off error.
		
\begin{figure}[h]
			\hfil 
			\includegraphics[width=8cm,height=6cm]{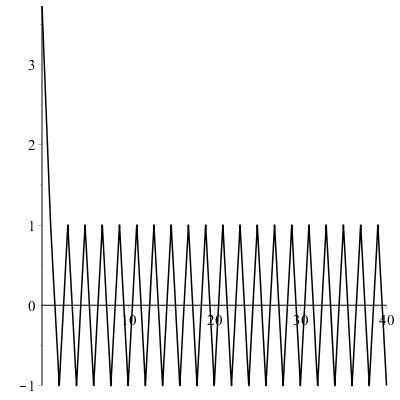}
			\caption{Choose $\theta_{0}=\dfrac{\pi}{12}$. Then we have $x_{1}=x_{3}$. The iteration starts to cycle from the second step.}
			\label{Hcycle}
		\end{figure}

		\paragraph{Case 2b} The initial value does not repeat. That is, the orbit is only \textsl{ultimately} periodic. Figure~\ref{Hcycle} is an example where $\theta_{0}=\pi/12$. By simple computation, it is easy to see that $x_{1}=\cot (\pi/4)= \sqrt{2}$, $x_{2}=\cot (3\pi/4)= -\sqrt{2}$, and $x_{3}=\cot (9\pi/4)= \sqrt{2}= x_{1}$. This is a period-$2$ cycle.   
		
		The map $\theta \to 3\theta$ mod $\pi$ is a Bernoulli shift. If we write the fraction $\dfrac{\theta}{\pi} = \dfrac{p}{q}$ in ternary, say $\dfrac{p}{q}=a_{0}.a_{1}a_{2}a_{3}a_{4}\dots$, then $\dfrac{3\theta}{\pi}$ moves the ternary point one place to the right, giving $a_{0}a_{1}.a_{2}a_{3}a_{4}\dots$. When the number is multiplied by $\pi$, the integer part makes no difference to the value of cotangent. So only the fractional part matters. 
		
		Look at all the examples again. For the one where the iteration blows up, the fraction is $1/9$ which is $0.01$ in ternary. This is a finite representation with two ternary places, while the sequence only exists for two steps. For the next two examples in case 2a, we notice that all the fractions can be represented by an infinite string of recurring digits in ternary. The fraction $1/8=0.\overline{01}$ has two digits in its repetend, while the iteration is a period-$2$ cycle satisfying $x_{0}=x_{2}$. Similarly, the fraction $1/7$ is $0.\overline{010212}$ in ternary, while the iteration is a period-$6$ cycle satisfying $x_{0}=x_{6}$. As for $1/12$ which is $0.0\overline{02}$ in ternary, there is a non-repeating digits right after the ternary point. So the example in case 2b starts to oscillate from the second step.

\begin{figure}[h]
			\hfil 
			\includegraphics[width=7cm,height=6cm]{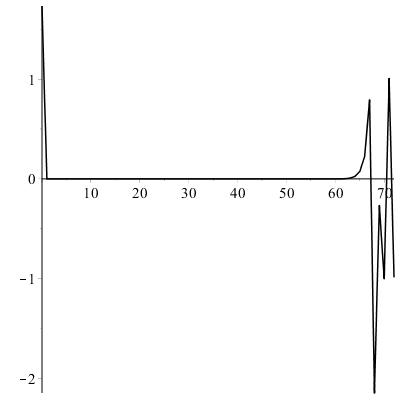}
			\caption{Take $\theta_{0}= \dfrac{\pi}{6}$ (its cotangent is $\sqrt{3}$). Then $x_{n}=0$ for all $n \geqslant 1$. The round-off error eventually pushes the sequence far away from zero.}
			\label{Hcycle3}
		\end{figure}

		\paragraph{Case 2c} A special case is the period-$1$ cycles, which means the iterations are actually convergent. If one of the steps returns zero, then iterates after that are zeros. This is obvious from the recursive formula $x_{n+1} = x_{n} \cdot \dfrac{x_{n}^{2}-3}{3x_{n}^{2}-1}$. However, the convergence is spurious since zero is not a root of $f(x)=x^2+1=0$. Notice that Halley's method is undefined since $f(0)=1$ and $f'(0)=0$, but $0-\dfrac{1}{0+(1/0)}$ could be interpreted as $0-\dfrac{1}{0+\infty}$. From the perspective of angles, if the cotangent of $\theta_{n}$ is zero, then $\theta_{n}$ and $\pi/2$ must differ by a multiple of $\pi$, and likewise $3\theta_{n}$ and $\pi/2$. Hence, the initial angle should be $\pm \dfrac{\pi}{2\cdot 3^{n}} (n \in \mathbb{Z})$. Choose $\theta_{0}= \dfrac{\pi}{6}=\dfrac{\pi}{2\cdot 3}$. Then $x_{0}=\sqrt{3}$ and $x_{n}=0$ for all $n \geqslant 1$. Figure~\ref{Hcycle3} shows the iteration numerically. The round-off error from the initial value grows with the computation, then eventually pushes the sequence far away from zero.
		
\begin{figure}[h]
			\hfil 
			\includegraphics[width=16cm,height=6cm]{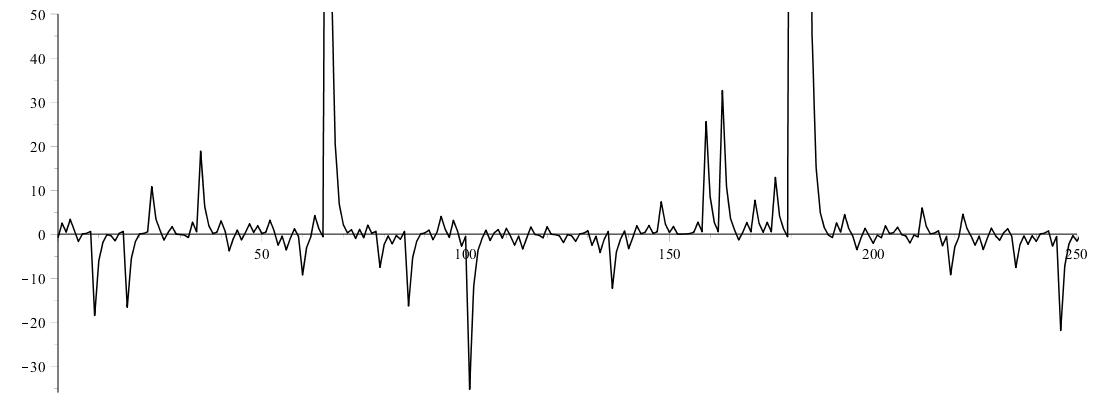}
			\caption{The iteration starting from $x_{0}= \cot\dfrac{\sqrt{2}\pi}{2}$ is not periodic. Two peaks are truncated to show more details around zero.}
			\label{Hchaos}
		\end{figure}

		\paragraph{Case 3} Considering the similarity between the formula for Newton's method and that for Halley's method, a good guess is that the iteration is not periodic if $\theta_{0}$ is an irrational multiple of $\pi$. Again, try $\theta_{0}= \dfrac{\sqrt{2}\pi}{2} \; (x_{0}=\cot \dfrac{\sqrt{2}\pi}{2})$. The first $250$ steps are shown in Figure~\ref{Hchaos}, which looks random but is in fact deterministic, corresponding to the ternary expansion of $1/\sqrt{2} = 0.2010021102221121\ldots$.
		
		It is easy to show that these sequences are aperiodic. Suppose that there exist two different terms which are equal to each other, say $x_{m} = x_{n}$. The two corresponding angles $\theta_{m}$ and $\theta_{n}$ differ by a multiple of $\pi$. Let $\theta_{0}=b\pi$. Then, $3^{m}b\pi = 3^{n}b\pi+k\pi \, (k \in \mathbb{Z})$, yielding $b=\dfrac{k}{3^{m}+3^{n}}$. Apparently, $b$ must be rational. Therefore, the sequence does not have repeating terms, if $\theta_{0}=b\pi$ for any irrational $b$.
		
\begin{figure}[h]
			\hfil
			\includegraphics[width=9cm,height=7.2cm]{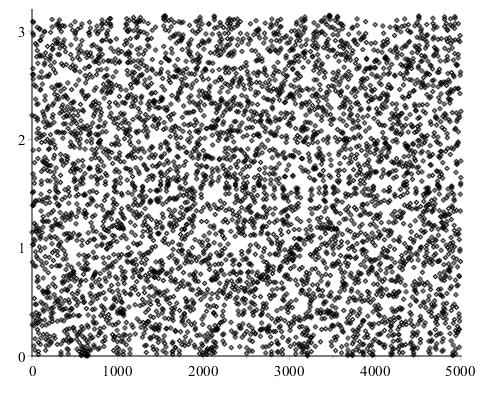}
			\caption{The sequence of angles modulo $\pi$ starting from $\theta_{0}=\dfrac{\sqrt{2} \pi}{2}$, looking apparently random.}
			\label{Hchaosangle}
		\end{figure}

		Considering the regular growth of the angles, we are more interested to the behaviour of the sequence $\{\theta_{n}\}$ modulo $\pi$. Use the same iteration above, the sequence of angles modulo $\pi$ is shown in Figure~\ref{Hchaosangle}.
		
		Naturally, the angles are bounded by $0$ and $\pi$. The sequence is non-periodic since we have proved that $\{x_{n}\}$ is not periodic. According to the expression $\theta_{n}=3^{n}\theta_{0}$, any little difference between the initial values will grow exponentially. In conclusion, this sequence is chaotic.

		\subsection{The secant method} \label{S::secant}
		
		Now we turn to the secant method. Different to Newton's method and Halley's method, this iteration is based on two previous steps. The recursive formula is 
		\begin{equation}
		x_{n+1} = x_{n}-\dfrac{f(x_{n})\left( x_{n}-x_{n-1} \right)}{f(x_{n})-f(x_{n-1})} = \dfrac{x_{n}x_{n-1}-1}{x_{n}+x_{n-1}},
		\end{equation}
		which looks similar to that for the cotangent of sums,
		\begin{equation}
		\cot{\left( \theta_{1}+\theta_{2} \right)}=\dfrac{\cot \theta_{1} \cot \theta_{2}-1}{\cot \theta_{1}+\cot \theta_{2}}.
		\end{equation}
		
		If $x_{n}=\cot{(\theta_{n})}$, $x_{n-1}=\cot{(\theta_{n-1})}$, then $x_{n+1}=\cot{(\theta_{n}+\theta_{n-1})}$. The list of angles is a general Fibonacci sequence. Therefore, if two initial points are given, namely $x_{0}=\cot{(\theta_{0})}$ and $x_{1}=\cot{(\theta_{1})}$, then
		\begin{equation}
		x_{n}=\cot{(F_{n-2}\theta_{0}+F_{n-1}\theta_{1})},
		\end{equation} 
		where $F_{n}$ denotes the $n$th term in the Fibonacci sequence.
		
		It is much more complicated for this formula to analyse the behaviour of iteration. But the results from above two methods suggest a way to try. 
		
		\paragraph{Guess 1} Given $\theta_{0}=a\pi$ and $\theta_{1}=b\pi$, if both $a$ and $b$ are rational numbers, then $x_{n}$ either diverges to infinity or eventually cycles. Here are two examples.
		\begin{itemize}
			\item Example 1: $\theta_{0}=\dfrac{\pi}{4}$ and $\theta_{1}=\dfrac{\pi}{2}$. 
			
			The sequence $\{\theta_{n}\}$ modulo $\pi$ is $\dfrac{\pi}{4}, \dfrac{\pi}{2}, \dfrac{3\pi}{4}, \dfrac{\pi}{4}, 0$. Hence, $\{x_{n}\}$ only exist for $n \leqslant 4$.
			
			\item Example 2: $\theta_{0}=\dfrac{\pi}{8}$ and $\theta_{1}=\dfrac{\pi}{2}$ 
			
			The initial values yield a sequence $\{ x_{n}\} $ of period $12$. 
		\end{itemize}  
		
		When the sequence goes to infinity, there exist $\theta_{N} \equiv 0 \mod \pi$. Let $\theta_{N-1} \equiv \theta \mod \pi$. Then $\theta_{N-2} \equiv \left( \pi-\theta\right) \mod \pi$. Iterating backwards to the initial points, we construct a new general Fibonacci sequence. Set $G_{-1}=\theta$ and $G_{-2}=\pi-\theta$. The Fibonacci recurrence can be written as $G_{-n}=G_{2-n}-G_{1-n}$. We obtain the general expression for the sequence,
		\begin{equation} 
		G_{-n}=F_{-n}\theta + F_{1-n}\pi=(-1)^{n+1}F_{n}\theta + (-1)^{n}F_{n-1}\pi,
		\end{equation}
		where $F_{n}$ denotes the $n$th term of the Fibonacci sequence. More details about the Fibonacci sequence identity can be found in Renault's work~(\cite{renault1996fibonacci}).
		
		Hence,
		\begin{align}
		G_{-N} &= (-1)^{N+1}F_{N}\theta + (-1)^{N}F_{N-1}\pi,   \notag \\
		G_{-(N-1)} &= (-1)^{N}F_{N-1}\theta + (-1)^{N-1}F_{N-2}\pi. \notag
		\end{align}
		Eliminating $\theta$, gives
		\begin{equation}  
		\frac{G_{-N}-(-1)^{N}F_{N-1}\pi}{(-1)^{N+1}F_{N}} = \frac{G_{-(N-1)}-(-1)^{N-1}F_{N-2}\pi}{(-1)^{N}F_{N-1}}.
		\end{equation}
		This equation can be simplified as below,
		\begin{equation}
		F_{N-1}G_{-N}+F_{N}G_{-(N-1)}=\pi,
		\end{equation}
		by using the identity $$F_{n}F_{n-2}-F_{n-1}^2=(-1)^{n-1}.$$
		
		Since $\theta_{0} \equiv G_{-N} \mod \pi$ and $\theta_{1} \equiv G_{-(N-1)} \mod \pi$, the initial angles must satisfy
		\begin{equation}
		F_{N-1}\theta_{0}+F_{N}\theta_{1} \equiv 0 \mod \pi,
		\end{equation}
		which is the condition for the iteration to blow up at step $N$. Otherwise, the iteration cycles. 
		
		\paragraph{Guess 2} Given $\theta_{0}=a\pi$ and $\theta_{1}=b\pi$, if either $a$ or $b$ is irrational, then $\{\theta_{n}\}$ modulo $\pi$ is aperiodic, so is $\{x_{n}\}$. The angles satisfy $\theta_{n+1}=\theta_{n}+\theta_{n-1}$.   
		
		\begin{figure}[!ht]
			\centering
			\subfigure[$\theta_{0}=\dfrac{\pi}{4}$, $\theta_{1}=\dfrac{\pi}{\sqrt{2}}$]{
				\label{subfig:Spoints:a} 
				\includegraphics[width=.475\textwidth,height=6cm]{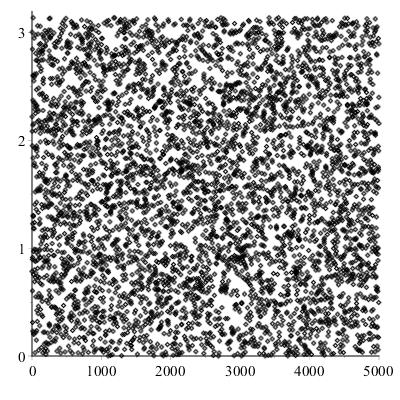}}\quad
			\subfigure[$\theta_{0}=\dfrac{\pi}{\sqrt{2}}$, $\theta_{1}=\dfrac{\pi}{4}$]{
				\label{subfig:Spints:b} 
				\includegraphics[width=.475\textwidth,height=6cm]{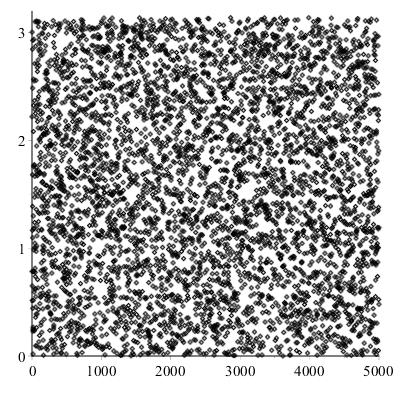}}\\
			\caption{The sequence of angles modulo $\pi$.}
			\label{fig:Spoints}
		\end{figure}
		
		Two examples are given in Figures~\ref{fig:Spoints}. The initial angles are $\theta_{0}=\pi/4$, $\theta_{1}=\pi/\sqrt{2}$ and $\theta_{0}=\pi/\sqrt{2}$, $\theta_{1}=\pi/4$, using the same angles but in different orders. Similar to the discussion about Halley's method Case~3, we can prove that this general Fibonacci sequence $\{\theta_{n}\}$ modulo $\pi$ is chaotic.
		
		The two guesses had been proved in Rhouma's work (see~\cite{rhouma2005fibonacci}, Theorem 1). In addition to their result, we have given the condition when the iteration blows up.

		\subsection{Householder methods}
		
		Generally, one can achieve arbitrary rate of convergence $k+1$ ($k \geqslant 1$), by the Householder method of order $k$~(\cite{householder1970numerical}), namely
		\begin{equation} \label{eqn::House}
		x_{n+1} = x_{n} + k \dfrac{(1/f)^{(k-1)}(x_{n})}{(1/f)^{(k)}(x_{n})}.
		\end{equation}
		Here, $F^{(n)}$ means the $n^{th}$ derivative of $F$. When $k=1$, this is just Newton's method since
		\begin{equation}
		\begin{aligned}
		x_{n+1} &= x_{n} + 1 \dfrac{(1/f)(x_{n})}{(1/f)^{(1)}(x_{n})} \\
		&= x_{n} + \dfrac{1}{f(x_{n})} \cdot \left( -\dfrac{f'(x_{n})}{f(x_{n})^{2}} \right)^{-1} \\
		&= x_{n} - \dfrac{f(x_{n})}{f'(x_{n})}.
		\end{aligned}
		\end{equation}
		When $k=2$, this is Halley's method since
		\begin{equation}
		\begin{aligned}
		x_{n+1} &= x_{n} + 2 \dfrac{(1/f)^{(1)}(x_{n})}{(1/f)^{(2)}(x_{n})} \\
		&= x_{n} + 2 \left( -\dfrac{f'(x_{n})}{f(x_{n})^{2}} \right) \cdot \left( \dfrac{-f''(x_{n})f(x_{n}) + 2f'(x_{n})^{2}}{f(x_{n})^{3}} \right)^{-1} \\
		&= x_{n} - \dfrac{2f'(x_{n})f(x_{n})}{2f'(x_{n})^{2} - f''(x_{n})f(x_{n})}.
		\end{aligned}
		\end{equation}
		When $k=3$, the rate of convergence is $4$. Iteration~(\ref{eqn::House}) becomes
		\begin{equation}
		\begin{aligned}
		x_{n+1} &= x_{n} + 3 \dfrac{(1/f)^{(2)}(x_{n})}{(1/f)^{(3)}(x_{n})} \\
		&= x_{n} + 3 \left( \dfrac{-f''(x_{n})f(x_{n}) + 2f'(x_{n})^{2}}{f(x_{n})^{3}} \right) \cdot \left( \dfrac{-f'''(x_{n})f(x_{n})^{2} + 6f''(x_{n})f'(x_{n})f(x_{n}) -6f'(x_{n})^{3}}{f(x_{n})^{4}} \right)^{-1} \\
		&= x_{n} - \dfrac{6f(x_{n})f'(x_{n})^{2} - 3f(x_{n})^{2}f''(x_{n})}{6f'(x_{n})^{3} - 6f''(x_{n})f'(x_{n})f(x_{n}) + f'''(x_{n})f(x_{n})^{2}}.
		\end{aligned}
		\end{equation}
		
		Substituting the function $f(x) = x^{2}+1$ and its derivatives, we obtain that
		\begin{equation}
		\begin{aligned}
		x_{n+1} &= x_{n} - \dfrac{(x_{n}^{2}+1)(3x_{n}^{2}-1)}{4x_{n}^{3}-4x_{n}} \\
		&= \dfrac{x_{n}^{4}-6x_{n}^{2}+1}{4x_{n}^{3}-4x_{n}},
		\end{aligned}
		\end{equation}
		which is similar to the cotangent identity,
		\begin{equation}
		\begin{aligned}
		\cot 4\theta &= \dfrac{\cot^{2} (2\theta) -1}{2 \cot (2\theta)} \\
		&= \dfrac{\cot^{4} \theta - 6\cot^{2} \theta +1}{4\cot^{3} \theta -4\cot \theta}.
		\end{aligned}
		\end{equation}
		Hence, if $x_{n}=\cot{\theta_{n}}$, then $x_{n+1}=\cot{(4\theta_{n})}$. Then, 
		\begin{equation} \label{eqnHouse}
		x_{n}=\cot{\left(4^{n}\theta_{0}\right)}, \quad \text{given} \quad x_{0}=\cot \theta_{0}.
		\end{equation} 
		This is equivalent to taking two Newton steps.
		
		
		\begin{theorem}\label{thm:House}
		The general Householder iteration of order~$k$ given in equation~\eqref{eqn::House} is solved by
		\begin{equation}
		    x_n = \cot\left( \left(k+1\right)^n \theta_0\right)
		\end{equation}
		where $\theta_0 \in (0,\pi)$ is determined by the initial condition $x_0 = \cot(\theta_0) \in \mathbb{R}$.
		\end{theorem}
		
\medskip\par\noindent		
\textbf{Proof.}
		\begin{equation}
		\dfrac{1}{f} = \dfrac{1}{x^{2}+1} = \dfrac{i/2}{x+i} - \dfrac{i/2}{x-i}.
		\end{equation} 
		Note that the $k^{th}$ derivative of $(x-a)^{-1}$ is 
		\begin{equation}
		\dfrac{(-1)^{k} \cdot k!}{(x-a)^{k+1}}.
		\end{equation}
		This trick for getting the symbolic $k^{th}$ derivative of a rational function is in~\cite{hardy2008course}, but is not generally taught in Calculus courses nowadays. Here,
		\begin{equation}
		\left( \dfrac{1}{f} \right)^{(k-1)} = \dfrac{(i/2)(-1)^{k-1}(k-1)!}{(x+i)^{k}} - \dfrac{(i/2)(-1)^{k-1}(k-1)!}{(x-i)^{k}},
		\end{equation}
		and similarly,
		\begin{equation}
		\left( \dfrac{1}{f} \right)^{(k)} = \dfrac{(i/2)(-1)^{k}k!}{(x+i)^{k+1}} - \dfrac{(i/2)(-1)^{k}k!}{(x-i)^{k+1}}\>.
		\end{equation}
		
		\textbf{Remark}.  Computer algebra systems have been able to do symbolic differentiation since the beginning.  Differentiation to a symbolic order is, of course, harder and came later.  All modern computer algebra systems are able to do this.  See for instance~\cite{GruntzKoepf::1995} or~\cite{Benghorbal::2002}. The result of the simple Maple command \verb`diff( 1/(x^2+1), x$n)` is equivalent to that above, although presented in a form that might be hard to read:
		\begin{equation}
		    \sum _{{\it \_alpha}={\it RootOf} \left( {{\it \_Z}}^{2}+1 \right) }-1
/2\,{\it \_alpha}\,{\mathrm{pochhammer}} \left( -n,n \right)  \left( x-{
\it \_alpha} \right) ^{-1-n}\>.
		\end{equation}
		
		Now back to the proof.  We consider the change of variable, 
		\begin{equation} \label{eqn::cot}
		x = \cot \theta = \dfrac{\cos \theta}{\sin \theta} = i \dfrac{e^{i \theta} + e^{-i \theta}}{e^{i \theta} - e^{-i \theta}}.
		\end{equation}
		Then,
		\begin{equation} \label{eqn::deviation}
		x+i = \dfrac{e^{i \theta}}{\sin \theta}, \quad \text{and} \quad x-i = \dfrac{e^{-i \theta}}{\sin \theta}.
		\end{equation}
		Thus, in the new variable,
		\begin{equation}
		k \dfrac{(1/f)^{(k-1)}(x_{n})}{(1/f)^{(k)}(x_{n})} = - \dfrac{\sin(k\theta_{n})}{\sin \theta_{n} \sin ((k+1) \theta_{n})}.
		\end{equation}
		Householder iteration then becomes,
		\begin{equation}
		\begin{aligned}
		\cot \theta_{n+1} &= \cot \theta_{n} - \dfrac{\sin(k\theta_{n})}{\sin \theta_{n} \sin ((k+1) \theta_{n})} \\
		&= \dfrac{\cos \theta_{n} \sin ((k+1) \theta_{n}) - \sin(k\theta_{n})}{\sin \theta_{n} \sin ((k+1) \theta_{n})} \\
		&= \cot ((k+1) \theta_{n}).
		\end{aligned}
		\end{equation}
		So we may take $\theta_{n+1} = (k+1) \theta_{n}$ mod $\pi$. This gives $\theta_{n} = (k+1)^{n} \theta_{0}$ mod $\pi$. For any real initial point $x_{0}$, there exist a unique $\theta_{0} \in (0, \pi)$ with $x_{0} = \cot \theta_{0}$. Then, as was to be proved
		\begin{equation}
		x_{n}=\cot{((k+1)^{n}\theta_{0})}\>.
		\end{equation} 

\par\medskip\noindent
\begin{flushright}
$\natural$
\end{flushright}
\par\medskip\noindent
		Similar to Newton iteration and Halley iteration, one can easily deduce the behaviour of a general Householder sequence $\{x_{n}\}$.
		\begin{enumerate}[(i)]
			\item If $\theta_{0}=\dfrac{M \pi}{(k+1)^{n}}$ for some $M \in \mathbb{Z}$, then $x_{n}=\cot{\left( M \pi \right)}$. The iteration blows up.
			\item If $\theta_{0}=\dfrac{p}{q} \pi$ for any fraction $\dfrac{p}{q}$ other than $\dfrac{M}{(k+1)^{n}} \left(M \in \mathbb{Z}\right)$, then the iteration eventually cycles. In addition, when $\theta_{0}=\dfrac{M \pi}{(k+1)^{n}-1}$ for some $M \in \mathbb{Z}$, we will see $x_{n}=x_{0}$. The iteration of period $n$ cycles from the start point.
			\item If $\theta_{0}=b \pi$ for some irrational number $b$, the iteration is not periodic (or convergent).
		\end{enumerate}

		Moreover, we can also prove the convergence of any complex sequences according to the deviations given by ~(\ref{eqn::deviation}). Denoting the complex initial point as $x_{0} = u+iv$, one can uniquely choose $\theta_{0} = \alpha + i\beta$, where $0 \leqslant \alpha \leqslant \pi$ and $x_{0} = \cot \theta_{0}$. Plugging into~(\ref{eqn::cot}) gives $v$ in terms of $\alpha$ and $\beta$,
		\begin{equation}
		v = \dfrac{(e^{-\beta} + e^{\beta})(e^{-\beta} - e^{\beta})}{(e^{-\beta} - e^{\beta})^{2} \cos^{2} \alpha + (e^{-\beta} + e^{\beta})^{2} \sin^{2} \alpha}.
		\end{equation}
		It is easy to verify that $v < 0$ if and only if $\beta>0$; $v > 0$ if and only if $\beta < 0$.
		
\begin{figure}[h]
			\hfil \includegraphics[width=6cm,height=6cm]{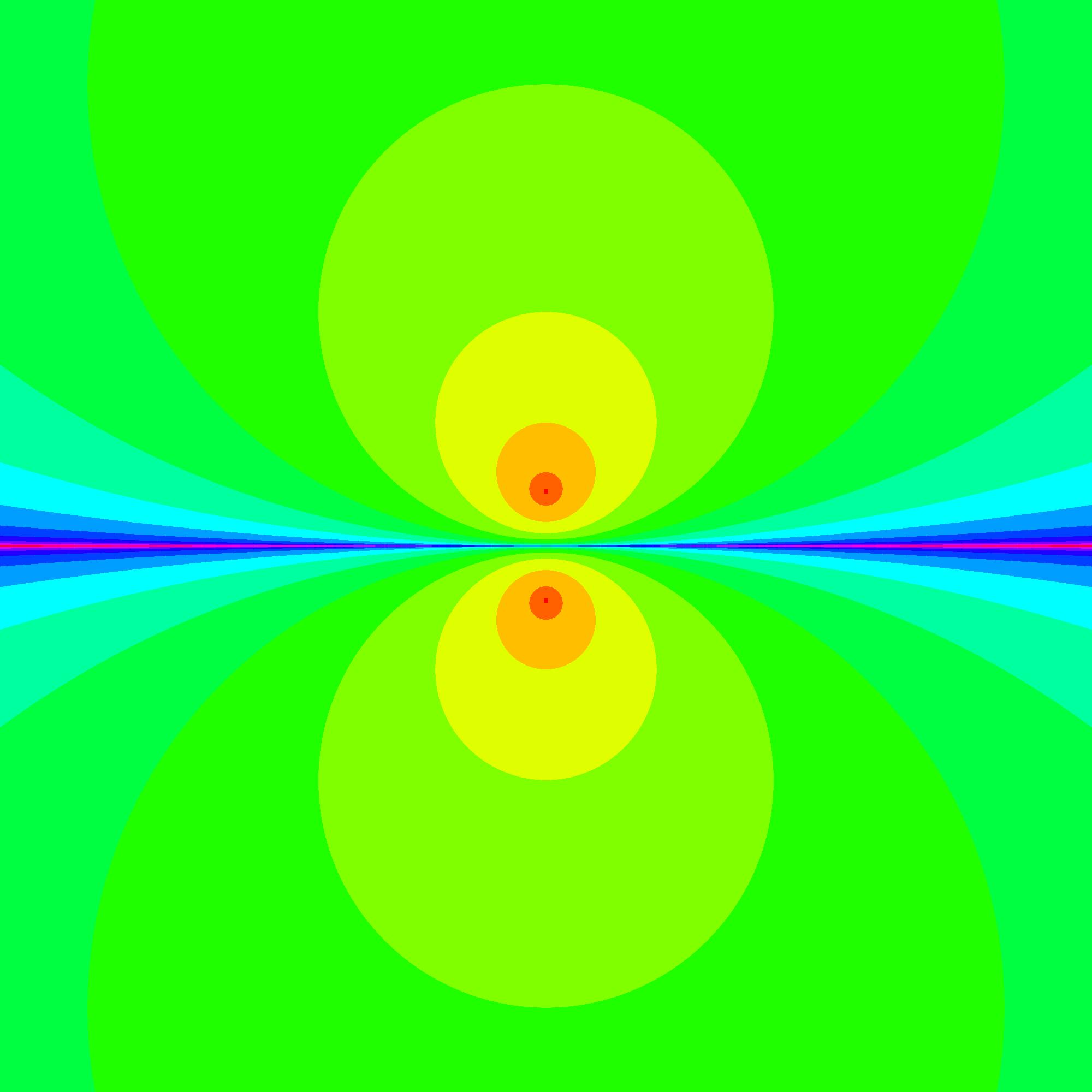}
			\caption{The Newton fractal of $f(x)=x^{2}+1$, where equation $f(x)=0$ has two roots $x=\pm i$. The basins of attraction are half-planes separated by the real axis.}
			\label{fig:fractal}
		\end{figure}

		The general Householder iteration of order $k$ gives $\theta_{n} = (k+1)^{n}(\alpha + i\beta)$. The deviations of $x_{n}$ from the roots $x = \pm i$ are,
		\begin{equation}
		x_{n} + i = \dfrac{e^{-\beta(k+1)^{n}} e^{i\alpha(k+1)^{n}}}{\sin ((k+1)^{n} (\alpha + i\beta))} \quad \text{and} \quad x_{n} - i = \dfrac{e^{\beta(k+1)^{n}} e^{-i\alpha(k+1)^{n}}}{\sin ((k+1)^{n} (\alpha + i\beta))}.
		\end{equation}  
		For any initial guess with $v > 0$, $x_{n}-i \to 0$ as $n \to \infty$, the sequence converges to the point $x=i$. For any initial guess with $v < 0$, $x_{n}+i \to 0$ as $n \to \infty$, the sequence converges to the point $x=-i$. The basins of attraction can be drawn as shown in Figure~\ref{fig:fractal}. All iterations starting from the upper semi-plane converge to $x=i$, while initial points in the lower half-plane lead towards $x=-i$. This diagram is well-known; see~\cite{kalantari2008polynomial} for beautiful generalizations.
		
		\section{Schr\"{o}der iterations are not so easy}
		
		In~\cite{petkovic2012schroder} we find a discussion showing that several classes of methods, including Householder's methods, are actually rediscoveries of Schr\"{o}der's second class of methods. By showing that Householder's methods give $x_{n}=\cot{(k+1)^{n} \theta_{0}}$ we have shown that all these methods (eight equivalent named classes of methods are given in~\cite{petkovic2012schroder}) give the same answers.
		
		Schr\"{o}der first class of methods is, however, not equivalent. We show below that Schr\"{o}der's first class of methods is unlikely to be explained by any equation similar to $x_{n} = S((k+1)^{n} \theta_{0})$ for any ``reasonable'' function $S$, at least for $k \geqslant 2$; for $k=1$, this method is also just Newton's method.
		
		\subsection{Reversion of series and Schr\"{o}der's first method}
		
		If $\Delta y$ has a Taylor series expansion in $\Delta x$, say $\Delta y = a_{1} \Delta x + a_{2} \Delta x^2 + a_{3} \Delta x^3 + \cdots$, then if $a_{1} \ne 0$ the expansion can be reversed (sometimes called ``reverted'') to get a series for $\Delta x$ in terms of $\Delta y$:
		\begin{equation}
		\Delta x = A_{1} \Delta y + A_{2} \Delta y^2 + A_{3} \Delta y^3 + \cdots.
		\end{equation}
		There are many treatments in the literature, and the idea goes back to Lagrange, and possibly to J. H. Lambert although his claim rests on his story that Acta Helvetica lost part of his manuscript; a beautiful algebraic exposition can be found in~\cite{henrici1974applied} , although Henrici there calls it the Lagrange-B\"{u}rman formula, whilst most authors just call it the Lagrange Inversion Formula.
		
		We do not need the full generality of these treatments, and can give instead the main idea of series reversion with the following simple computation: we put the known series for $\Delta y$ in terms of $\Delta x$ into the reverted series, and equate powers of $\Delta x$. (It works just as well if we put the reverted series into the original.)
		\begin{equation}
		\begin{aligned}
		\Delta x &= A_{1} (a_{1} \Delta x + a_{2} \Delta x^2 + \cdots) + A_{2} (a_{1}^2 \Delta x^2 + \cdots) + \cdots \\
		&= A_{1}a_{1}\Delta x + (A_{1}a_{2} + A_{2}a_{1}^2)\Delta x^2 + \cdots
		\end{aligned}
		\end{equation} 
		Obviously, 
		$$A_{1} = \dfrac{1}{a_{1}} \quad \text{and} \quad A_{2} = - \dfrac{A_{1}a_{2}}{a_{1}^2} = -\dfrac{a_{2}}{a_{1}^3}.$$
		
		One can carry this argument out to any desired order, and indeed the first few results are even tabulated in \cite{abramowitz1967handbook} (page 16).
		Nowadays one prefers to use computer algebra, and in Maple the simplest thing is to use the \texttt{solve} command on a \texttt{series}.
		For instance, if the variable \texttt{Order} is set to $4$ and the variable $Y$ contains a \texttt{series} 
		\begin{equation}
		    Y = \eta+a_{{1}} \left( x-\xi \right) +a_{{2}} \left( x-\xi \right) ^{2}+
a_{{3}} \left( x-\xi \right) ^{3}+O \left(  \left( x-\xi \right) ^{4}
 \right)\>,
		\end{equation}
and one issues the command \verb`solve(y=Y,x)`, one gets 
\begin{equation}
\xi-{a_{{1}}}^{-1} \left( \eta-y \right) -{\frac {a_{{2}}}{{a_{{1}}}^
{3}}} \left( \eta-y \right) ^{2}+{\frac {a_{{1}}a_{{3}}-2\,{a_{{2}}}^{
2}}{{a_{{1}}}^{5}}} \left( \eta-y \right) ^{3}+O \left(  \left( \eta-y
 \right) ^{4} \right)
\end{equation}
for $x$.  This is correct, although it would have been nice to have an expansion in $y-\eta$ automatically (one can get this by calling \texttt{series} on the result, and this fixes all the signs).

		For this specific application, we seek a zero of $y=f(x)$. We expand about our guess $x_{n}$:
		\begin{equation}
		y=f(x_{n}) + f'(x_{n})(x-x_{n}) + \dfrac{f''(x_{n})}{2}(x-x_{n})^2 + \cdots.
		\end{equation}
		Put $\Delta y = y-f(x_{n})$ and $\Delta x = x-x_{n}$. Then,
		\begin{equation}
		\Delta y = f'(x_{n})\Delta x + \dfrac{f''(x_{n})}{2} \Delta x^2 + \cdots,
		\end{equation}
		and $a_{1} = f'(x_{n})$, $a_{2} = f''(x_{n})/2$, etc. Reversion gives
		\begin{equation}
		\Delta x = A_{1} \Delta y + A_{2} \Delta y^2 + \cdots,
		\end{equation}
		where $A_{1} = 1/f'(x_{n})$, $A_{2} = -f''(x_{n})/(2f'(x_{n})^3)$, etc. Now, we are looking for $x$ so that $y=0$; then,
		\begin{align}
		\Delta y &= 0 - f(x_{n}) = -f(x_{n}), \quad \text{and}\\
		\Delta x &= \dfrac{1}{f'(x_{n})} (-f(x_{n})) - \dfrac{f''(x_{n})}{2f'(x_{n})^3} (-f(x_{n}))^2 + \cdots.
		\end{align}
		Truncations of these various reversions give Schr\"{o}der's first class of iterations: $\Delta x = x_{n+1} - x_{n}$ and Schr\"{o}der's third order method is simply
		\begin{equation}
		x_{n+1} - x_{n} = -\dfrac{f(x_{n})}{f'(x_{n})} - \dfrac{f''(x_{n}) f(x_{n})^2}{2f'(x_{n})^3}.
		\end{equation}
		For $f(x) = x^2+1$ this gives, after some algebra,
		\begin{equation}
		\begin{aligned}
		x_{n+1} &= G(x_{n}) \\
		&= x_{n} - \dfrac{5x_{n}^4 + 6x_{n}^2 + 1}{8x_{n}^3} \\
		&= \dfrac{3x_{n}^4 - 6x_{n}^2 - 1}{8x_{n}^3}. 
		\end{aligned}
		\end{equation}
		We will use this to show that Schr\"{o}der's third order method gives an iteration too complicated to explain with $x_{n} = S(3^n \theta_{0})$ for any reasonable function $S$. 
		
\begin{figure}[h]
			\hfil \includegraphics[width=6cm,height=6cm]{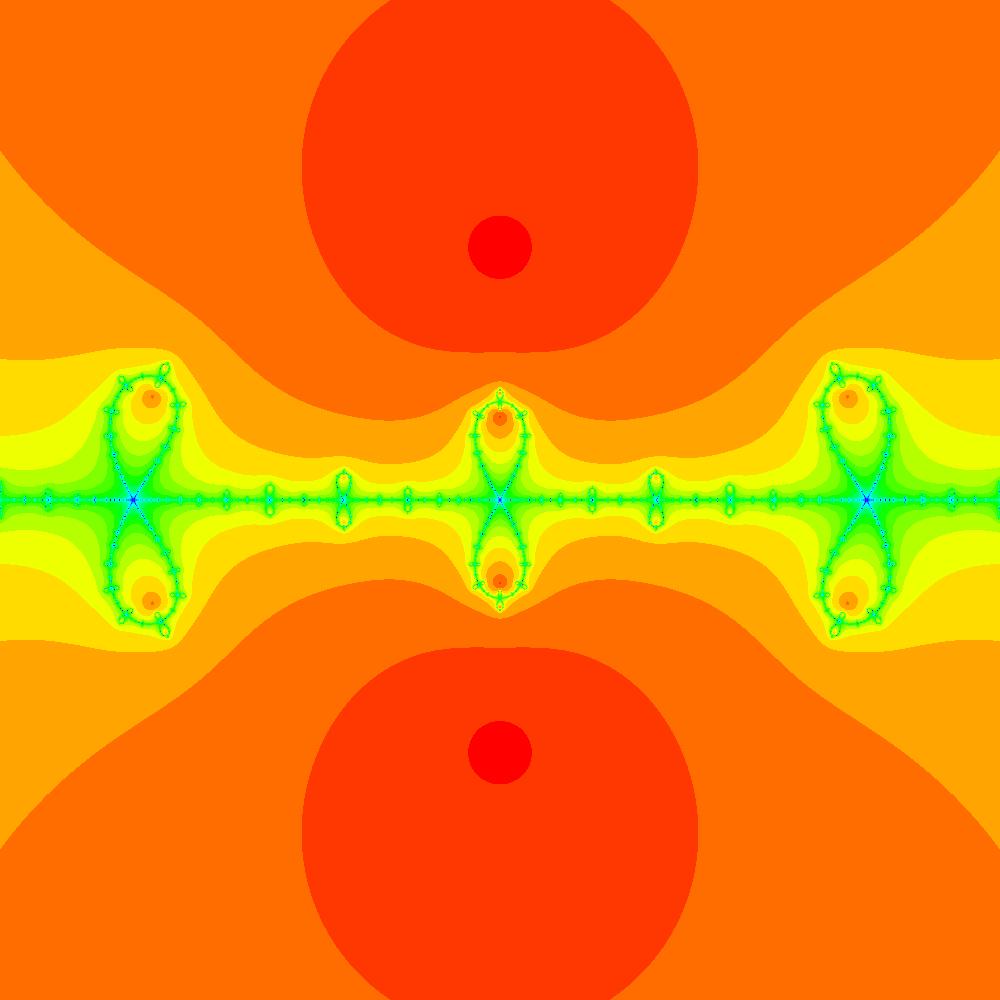}
			\caption{The basins of attraction of the roots of $f(x)=x^{2}+1=0$ by using Schr\"{o}der's third order iteration.}
			\label{fig:Schroder}
		\end{figure}

		The first thing we do is derive an equivalent function $F(x)$ for which the iteration above, $x_{n+1} = G(x_{n})$, is Newton's iteration. The function $F(x)$ satisfies
		\begin{align}
		\dfrac{F(x)}{F'(x)} &= \dfrac{5x_{n}^4 + 6x_{n}^2 + 1}{8x_{n}^3} \quad \text{or}  \\
		\dfrac{F'(x)}{F(x)} &= \dfrac{8x_{n}^3}{(5x^2+1) (x^2+1)} = -\dfrac{1}{5} \dfrac{10x}{5x^2+1} + \dfrac{2x}{x^2+1}. 
		\end{align}
		Integrating both sides yields
		\begin{equation}
		\ln F(x) = -\dfrac{1}{5} \ln (5x^2+1) + \ln (x^2+1).
		\end{equation}
		Constants of integration are immaterial here. Thus,
		\begin{equation}
		F(x) = (x^2+1)(5x^2+1)^{-1/5}.
		\end{equation}
		That is, Schr\"{o}der's third order iteration on $x^2+1$ is exactly Newton iteration on $(x^2+1)(5x^2+1)^{-1/5}$. This allows us to use the computer algebra system Maple to (quickly) draw the basins of attraction of the roots at $x=\pm i$. See Figure~\ref{fig:Schroder}.
		
		Notice that the iteration has two spurious fixed points: $x=G(x)$ implies $(5x^2+1)(x^2+1)=0$ which is possible not only when $x=\pm i$ but also when $x=\pm i/\sqrt{5}$. For the latter, $G'(\pm i/\sqrt{5}) = -3/2$ which is larger than $1$ in magnitude so these fixed points are repelling. Hence, the basins in Figure~\ref{fig:Schroder} have (in our opinion, beautiful) fractal boundaries. 
		
		Let us now consider what this means. If there were a simple function $S$ such that $x_{n} = S(3^n \theta_{0})$, then for some $\theta_{0}$, namely those with $x_{0} = S(\theta_{0})$ inside the basin of attraction of $i$, we would have $S(3^n \theta_{0}) \to i$; likewise with some other $\theta_{0}$, namely those with $x_{0} = S(\theta_{0})$ inside the basin of attraction of $-i$, we would have $S(3^n \theta_{0}) \to -i$. Thus, the function $S(\theta)$ would inherently contain information about the fractal boundary pictured in Figure~\ref{fig:Schroder}.
		
		For us to have a formula with $x_{n} = S(\theta)$ and $x_{n+1} = S(3\theta)$ we must have 
		\begin{equation} \label{eqn::S}
		S(3\theta) = G(S(\theta)),
		\end{equation}
		a functional equation for the unknown $S(\theta)$. Moreover, if $S(\theta) = \pm i$, $S(3\theta)$ must also be $\pm i$, and similarly if $S(\theta) = \pm i/\sqrt{5}$, then $S(3\theta) = \pm i/\sqrt{5}$ also. We have been unable to solve this functional equation. It is certainly true that $S(\theta) = \cot \theta$ does not solve it. If we look for functions $S(\theta)$ with algebraic singularities at $\theta = 0$, $S(\theta) = c \cdot \theta^{-\alpha} + o(\theta^{-\alpha})$ as $\theta \to 0$ for some $\alpha > 0$, then condition~(\ref{eqn::S}) requires
		\begin{equation}
		c(3\theta)^{-\alpha} + o(\theta^{-\alpha}) = \dfrac{3}{8} c \cdot \theta^{-\alpha} + o(\theta^{-\alpha})
		\end{equation}
		which can only be true if $3^{-\alpha} = 3/8$ or $$\alpha = \dfrac{\ln (8/3)}{\ln 3} \approx 0.892789.$$ This rules out many simple elementary functions already.
		
		Similarly if $S(\theta)$ has a logarithmic behaviour, $S(\theta) \sim \alpha \ln \theta + o(\ln \theta)$ perhaps, then
		\begin{equation}
		\alpha \ln \theta + \alpha \ln 3 + \cdots = \dfrac{3}{8} \alpha \ln \theta + \cdots
		\end{equation}
		which is impossible unless $\alpha = 0$.
		
		These computations do not (as far as we know!) prove that such an $S(\theta)$ is not elementary; but they suggest that Schr\"{o}der iterations are more difficult to analyze for this problem than Householder iterations. We conclude that the behaviour on the real axis is unlikely to be described simply. We would be interested in any clarification that might be provided by expert readers. Can equation~\eqref{eqn::S} be solved by an elementary $S$?

		\section{Discussion}
		
		Iteration of simple functions can produce complex behaviour. For instance, the well-studied quadratic iteration $z_{n+1}=az_{n}\left( 1-z_{n} \right) $ leads to chaos~\cite{gleick2011chaos}. We believe this present paper will help to understand the dynamic behaviour of chaos in another way. Besides, when using these classical numerical methods, such as Newton's iteration, Halley's iteration and the secant iteration, one needs to be aware of that these methods can fail.
		
		Another fact of note is that all the analytical expressions are related to the rates of convergence. The formulas for Newton's method, Halley's method and the secant method are 
		\begin{equation} \notag
		\begin{aligned}
		x_{n} &= \cot{(2^{n}\theta_{0})}, \\ 
		x_{n} &= \cot{(3^{n}\theta_{0})}, \\ 
		x_{n} &= \cot{(F_{n-2}\theta_{0}+F_{n-1}\theta_{1})}, \quad \text{where} \quad F_{n}=\dfrac{1}{\sqrt{5}} \left[ \left(\dfrac{1+\sqrt{5}}{2}\right)^n - \left(\dfrac{1-\sqrt{5}}{2}\right)^n \right],
		\end{aligned}
		\end{equation}
		given $x_{0}=\cot{\theta_{0}}$. And their rates of convergence are $2$, $3$ and $\dfrac{1+\sqrt{5}}{2}$, respectively. We also proved that the iteration for Householder's method with rate of convergence $k$ is $x_{n}=\cot{(k^{n}\theta_{0})}$. However, neither Schr\"{o}der's first method nor the basic sequence of Kalantari give cotangent formulas that we could find.
		
		On the other hand, any one-step iteration is Newton's method~(\cite{bateman1938haley}). In the last section, we used this idea to draw the basins of attraction for Schr\"{o}der's iteration. This can be extended to any scalar iteration,
		\begin{equation}
		x_{n+1} = H(x_{n}),
		\end{equation} 
		which is equivalent to Newton's iteration for function $h(x)$ if 
		\begin{equation}
		x - \dfrac{h(x)}{h'(x)} = H(x)
		\end{equation}
		for all $x$. This is a differential equation for $h(x)$, given $H(x)$; moreover, it is separable:
		\begin{equation}
		x - H(x) = \dfrac{h(x)}{h'(x)} 
		\end{equation}
		or
		\begin{equation}
		\dfrac{1}{x - H(x)} = \dfrac{h'(x)}{h(x)}, \quad \text{if} \quad x - H(x) \ne 0.
		\end{equation}
		Integrating both sides yields
		\begin{equation}
		h(x) = h(x_{0}) e^{\displaystyle{\int_{x_{0}}^{x}} \frac{d\xi}{\xi - H(\xi)}}.
		\end{equation}
		
		
		This fact, that any one-step iteration is equivalent to a Newton iteration for some other scalar function, is frequently rediscovered. The earliest reference we know for this is~\cite{bateman1938haley}. The most recent reference connecting iterations to Newton's method that we know is~\cite{tapia2018inverse}, where the authors carefully extend this idea to systems.
		
		Simulating mathematical dynamical systems in floating-point arithmetic can give surprising differences to what is expected. In this paper we have given some new mathematical analyses of dynamical systems that arise when using root-finding methods on a simple equation. Similar behaviour can occur for more complicated equations. We have also confirmed by example that floating-point arithmetic can alter the predicted behaviour. Of course, owing to the exponential sensitivity of chaotic systems, this is to be expected.
		
		We have not analyzed in detail the effect of floating-point arithmetic on these examples, as was done in~\cite{corless1992continued} for the Gauss map; we believe that this could be done, and a similar ``shadowing'' result proved---essentially constructing the ternary or $(k+1)$-ary expansion of $\theta_{0}/\pi$ retrospectively from the computed orbit---but we have not done so. A more intriguing question that remains is just how representative of true reality are these computed shadows? We leave that question for a future investigation. 
		
		Gilbert Strang's delightful article~\cite{strang1991chaotic} is very informative about Newton's method, chaos, and the power of exact solutions. This present paper only pushes those insights a little further. It is not really surprising that Schr\"{o}der's first (third order) method is not as simple as Newton's method; it is quite surprising that Halley's method, Householder's methods, and the secant method are in fact just as simply explained. We hope, however, that you (the readers) have gained some appreciation of the scope of research into root-finding methods, and the power of computer algebra systems to do so, even with this simple example.
		
	 \section*{Acknowledgement} 
			
			We thank David W.~Linder at Maplesoft for help with the \textsl{Fractals} package in Maple. This work was supported by the Natural Science and Engineering Research Council of Canada.  Support from the Rotman Institute for Philosophy and the School of Mathematical and Statistical Sciences at Western, and from the Ontario Research Centre for Computer Algebra (ORCCA) is gratefully acknowledged.

		\vfill\eject
	
\end{document}